\documentclass[11pt]{amsart}

\usepackage{enumerate}
\usepackage{amssymb}
\usepackage{graphics}
\usepackage{graphicx}
\usepackage{amscd}
\usepackage{amsmath, epsfig,  psfrag}
\usepackage{pinlabel}

\newtheorem{theorem}{Theorem}
\newtheorem{lemma}[theorem]{Lemma}
\newtheorem{prop}[theorem]{Proposition}

\newtheorem{cor}[theorem]{Corollary}

\newtheorem{convention}[theorem]{Convention}

\theoremstyle{definition}
\newtheorem{remark}[theorem]{Remark}
\newtheorem{example}[theorem]{Example}

\numberwithin{theorem}{section}

\newcommand{\on}{\operatorname}
\newcommand{\DD}{{\bf D}}

\renewcommand{\d}{\partial}
\newcommand{\Map}{{\on{Map} \DD_n}}

\newcommand{\ZZ}{\mathbb{Z}}

\newcommand{\tb}{\on{tb}}

\begin{document}

\author{Olga Plamenevskaya}
\address{Department of Mathematics, Stony Brook University, Stony Brook, NY 11794}
\email{olga@math.sunysb.edu}

\thanks{The first author is partially supported by 
NSF grant DMS-0805836.}

\author{Jeremy Van Horn-Morris}
\address{American Institute of Mathematics, Palo Alto, CA 94306}
\email{jvanhorn@aimath.org}

\title[Planar open books, monodromy factorizations, and symplectic fillings]{Planar open books, monodromy factorizations, \\ and symplectic fillings}

\begin{abstract} 

We study fillings of contact structures supported by planar open books
by analyzing positive factorizations of their monodromy. Our method is based on 
Wendl's theorem on symplectic fillings of planar open books.
We prove that every virtually overtwisted contact structure on $L(p,1)$ has a unique filling,  
and describe fillable and non-fillable tight contact structures on certain Seifert fibered spaces. 

\end{abstract}
\maketitle

\section{Introduction}  

By Giroux's theorem \cite{Gi}, a  contact 3-manifold $(Y, \xi)$ is Stein fillable if and only if 
it is compatible with an open book $(S, \phi)$ whose monodromy  $\phi$ can be represented as a  product of positive Dehn twists. Given a factorization  of the monodromy into a product $\phi= D_{\alpha_1} \dots D_{\alpha_k}$ of positive Dehn twists around homologically non-trivial curves $\alpha_1, \dots, \alpha_k$, we can construct a Stein filling as an allowable Lefschetz fibration over $D^2$ with fiber $S$, with vanishing cycles corresponding to  $\alpha_1, \dots, \alpha_k$. (We say that a Lefschetz fibration is allowable if all vanishing cycles are homologically non-trivial in their fibers.) 
 Conversely, if $X$ is a Stein manifold whose boundary is  $(Y, \xi)$, then $X$ has a structure of an allowable  Lefschetz fibration \cite{AO}. The boundary $\d X=Y$ has an open book decomposition whose monodromy is a product of positive Dehn twists around curves corresponding to the vanishing cycles; this open book is compatible with the contact structure $\xi$ \cite{Pl}. Thus, Stein fillings of $(Y, \xi)$ correspond to positive factorizations of monodromies of compatible open books. However, to detect non-fillability or to classify all Stein fillings, one would have to consider {\em all} possible open books compatible with $\xi$. 

The situation is much simpler for contact structures compatible with planar open books, thanks to the following recent result  
of Chris Wendl.
\begin{theorem}[Wendl \cite{We}] \label{wendl} Suppose that  $(Y,\xi)$ admits a planar open book decomposition. 
Then every strong symplectic filling $(X, \omega)$ of $(Y,\xi)$ is symplectic deformation equivalent 
to a blow-up of an allowable Lefschetz fibration compatible with the given open book for $(Y, \xi)$. 
\end{theorem} 
In particular, Wendl's theorem implies that every Stein filling of $(Y, \xi)$ is diffeomorphic (and even symplectic deformation 
equivalent) to an allowable Lefschetz fibration compatible with the given planar open book; to classify fillings or to prove non-fillability, 
it suffices to study positive factorizations of the {\em given} monodromy. Even so, enumerating positive factorizations for a given 
element of the mapping class is in general a very hard question. However, we are able to analyze certain simple monodromies by 
means of elementary calculations in the abelianization of the mapping class group of the planar surface. In this paper, 
we give two applications: first, we complete  the classification of fillings for tight  lens spaces $L(p,1)$,  second,  
we study the fillability question for certain tight Seifert fibered spaces. 

\begin{theorem} Every virtually overtwisted contact structure on $L(p,1)$ has a unique Stein filling (up to symplectic deformation), which is also its unique weak symplectic filling (up to symplectic deformation and blow-up).   \label{main} 
\end{theorem}     
   
%
%
%
%
%
%

\begin{cor} For $p\neq 4$, every tight contact structure on  $L(p,1)$ has a unique Stein filling (up to symplectic deformation), which is also its unique weak symplectic filling (up to symplectic deformation and blow-up).
\end{cor}

The above corollary combines Theorem \ref{main} together with earlier results, giving a complete description of fillings for tight $L(p,1)$.
(Note that fillings for $L(4,1)$ are also understood due to \cite{McD} and our theorem.)
Recall that Eliashberg established uniqueness of a symplectic filling (up to deformation and blow-up) for the standard contact structure on $S^3$ \cite{El}.  
McDuff  proved that standard contact structures on $L(p,1)$ all have unique filling except for $L(4,1)$, which has two fillings, up to blow-up and diffeomorphism \cite{McD}. Hind  showed that in McDuff's theorem, 
the Stein filling of $(L(p,1), \xi_{std})$ is 
 in fact unique up to Stein homotopy \cite{Hi}.  Lisca extended these results to obtain a classification (again up to blow-up and diffeomorphism) of symplectic fillings for arbitrary lens spaces $L(p,q)$ equipped with standard contact 
structures \cite{Li}. (By the standard contact structure on a lens space we mean the quotient of $(S^3, \xi_{std})$ by the action of the cyclic group. All of the above  results for the standard contact structure obviously extend to its conjugate, thus covering the two universally tight contact structures on $L(p,q)$.) Theorem \ref{main} extends classification of fillings in another direction, 
by a different technique. Our technique also allows to reprove uniqueness of symplectic fillings of $(L(p,1), \xi_{std})$ for
$p\neq 4$.  
\begin{figure}[htb] 
	\labellist
	\small\hair 2pt
	\pinlabel {$-1$} at 137 68
	\pinlabel {$-\frac{1}{r_1}$} at 37 37
	\pinlabel {$-\frac{1}{r_2}$} at 80 30
	\pinlabel {$-\frac{1}{r_3}$} at 122 36
	\pinlabel {$-1$} at 366 68
	\pinlabel {$-a_1$} at 265 64
	\pinlabel {$-a_2$} at 261 32
	\pinlabel {$-a_{n_1}$} at 252 2

	\pinlabel {$-b_1$} at 305 62
	\pinlabel {$-b_2$} at 309 29
	\pinlabel {$-b_{n_2}$} at 314 -1

	\pinlabel {$-c_1$} at 329 68
	\pinlabel {$-c_2$} at 337 32
	\pinlabel {$-c_{n_3}$} at  344 9

\endlabellist
\centering	
	
\includegraphics[scale=1.0]{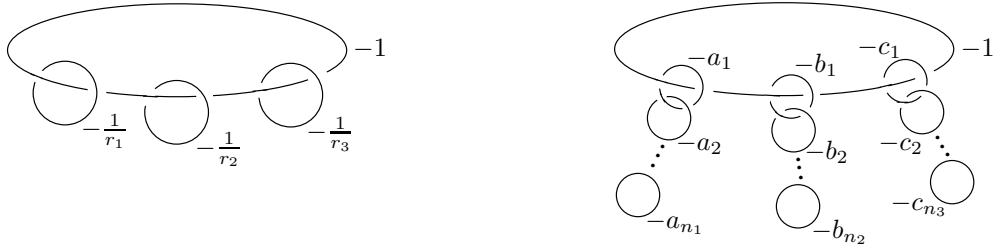}

\caption{Seifert fibered space  $M(-1; r_1, r_2, r_3)$.}\label{seifert} 
\end{figure}  

Our second application concerns fillability of contact structures on Seifert fibered spaces  $M(-1; r_1, r_2, r_3)$. (We use this notation for 
the space given by the surgery diagram of Figure \ref{seifert}; here and throughout 
the paper, $r_1, r_2, r_3$ are rational numbers between 0 and 1.) Tight contact structures on such manifolds were studied by 
Ghiggini--Lisca--Stipsicz \cite{GLS} and Lisca--Stipsicz \cite{LS, LS2}; when $r_1, r_2\geq\frac12$, a complete classification 
of tight contact structures on  $M(-1; r_1, r_2, r_3)$ was obtained in \cite{GLS}  (in particular, each of these spaces is known to carry
a tight contact structure). Tightness of some of these contact structures was 
established by means of the Heegaard Floer theory; it was shown in \cite{GLS} that 
one of the tight structures on $M(-1; \frac{1}2, \frac12, \frac1p)$ is non-fillable. (Recall that, in contrast, all tight contact structures 
on $M(0; r_1, r_2, r_3)$  are fillable, see \cite{GLS2}, cf \cite{Wu}).   

It is interesting to determine which of the manifolds $M(-1; r_1, r_2, r_3)$ carry tight, non-fillable contact structures.  Wendl's theorem 
provides a good tool for this investigation, because all tight contact structures on $M(-1; r_1, r_2, r_3)$ admit planar open books, 
at least in the case $r_1, r_2\geq\frac12$. It turns out that fillability depends, in a rather subtle way, on the arithmetics 
of the continued fraction expansions of $r_1, r_2, r_3$. Let 
$$
\displaystyle{-\frac1{r_1}= [a_1, a_2, \dots, a_{n_1}]},  \quad
\displaystyle{-\frac1{r_2}= [b_1, b_2, \dots, b_{n_2}]}, \quad
\displaystyle{-\frac1{r_3}= [c_1, c_2, \dots, c_{n_3}]},
$$ where we adopt the notation
$$
\displaystyle{[x_1, x_2, \dots, x_{n}] = - x_1 - \frac{1}{\displaystyle{-x_2 -\displaystyle{\frac{1}{\ddots- \displaystyle{ \frac{1}{-x_{n}}}    }    }}}}
,  \qquad x_i \in \ZZ, \ \ x_i \geq 2.
$$

\begin{theorem} \label{th2}  Suppose that $r_1, r_2\geq\frac12$. Let $k_1$, $k_2$ be such that $a_1=a_2=\dots=a_{k_1}=2$,  $b_1=b_2=\dots=b_{k_2}=2$, and $a_{k_1+1} \geq 3$, $b_{k_2+1} \geq 3$ (if $k_1<n_1$, resp. $k_2<n_2$). 
Then the space  $M(-1; r_1, r_2, r_3)$ carries tight, symplectically non-fillable contact 
structures if $c_1-1> \max(k_1,k_2)$; otherwise all tight structures on $M(-1; r_1, r_2, r_3)$  are Stein fillable. 
\end{theorem}

It is instructive to compare our method to other proofs of non-fillability. Previous techniques (see e.g. \cite{Li1, GLS}) are all based on some version 
of gauge theory and the Donaldson thorem, as follows. First one shows, using Seiberg-Witten or Heegaard Floer theory, 
that all symplectic fillings of a given contact structure must be negative-definite; then the filling is completed to a closed negative-definite 
4-manifold, so that by  the Donaldson theorem the intersection form of the filling embeds into a standard diagonalizable form over the integers. 
Finally, existence of such an embedding must be ruled out; often, this is possible because the first Chern class of a Stein filling can be understood 
in terms of the contact structure. We note that this method is in principle applicable to planar open books, because by \cite{Et2, OSS} 
every symplectic filling of a contact structure $\xi$ supported by a planar open book must be negative-definite, and have the intersection form 
that embeds  into a standard diagonalizable form over the integers; moreover, if $c_1(\xi)=0$, then $c_1$ of every Stein filling of $\xi$ must be 0.
However, the analysis of possible embeddings of the intersection form can typically be done only in a very limited number of cases. 
Our analysis of the monodromy factorizations, at least on the level of the abelianization, appears be a lot easier to do, and works in many situations.

\subsection*{Acknowledgements} We are grateful to John Etnyre for some helpful conversations,  to 
Dan Margalit and Chris Wendl for  helpful correspondence, and to Paolo Lisca and Andr\'as Stipsicz
for their comments on a preliminary version of this paper.
\section{Lens spaces $L(p,1)$ and their fillings}

By Honda's classification of tight contact structures on lens spaces \cite{Ho}, all tight contact structures on $L(p,1)$ arise 
as surgeries on Legendrian unknot with $tb=-p+1$, i.e. the standard Legendrian unknot stabilized $p-2$ times. More precisely, 
there are $p-2$ tight contact structures $\xi_1, \xi_2, \dots \xi_{p-2}$ on $L(p,1)$, distinct up to isotopy; the contact structure 
$\xi_k$ is the result of Legendrian surgery on the stabilized unknot with $k$ cusps on the left and $p-k$ cusps on the right.
Legendrian surgery yields a Stein filling for $\xi_k$ which is diffeomorphic to the union of $D^4$ and the 2-handle corresponding to the surgery.

Our job is to prove that every Stein filling for  a virtually overtwisted contact structure on  $(L(p,1), \xi)$ is  
diffeomorphic to the one described above, and that the Stein structure is unique, at least up to symplectic deformation. 
In our notation, the virtually overtwisted contact structures are $\xi_2, \dots, \xi_{p-3}$;  
indeed, by \cite{Go},  the two universally tight contact structures are given by surgeries on the unknot that has all of the stabilizations on the right or on the left, and the virtually overtwisted contact structures correspond to unknots stabilized on both sides.

\begin{figure}[htb] 
\includegraphics[scale=0.8]{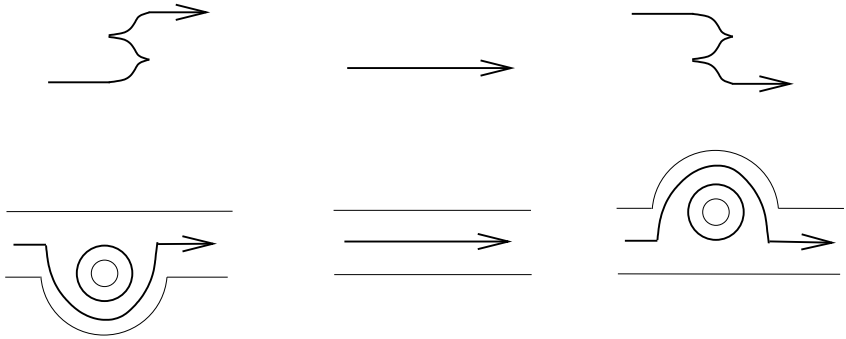}
\caption{Placing a stabilized knot on a page of an open book.}\label{stab} 
\end{figure}

We construct planar open books for  $(L(p,1), \xi_k)$ as follows.
The standard  Legendrian unknot with $\tb=-1$ can be represented as the core circle of a page of the open book 
decomposition of $(S^3, \xi_{std})$ with annular pages; to place the stabilized unknot on a page of an open book for  $(S^3, \xi_{std})$, we stabilize 
the open book and modify the knot as described in \cite{Etn}, see Figure \ref{stab}.

\begin{figure}[htb] 
	\labellist
	\small\hair 2pt
	\pinlabel {\mbox{$k$ cusps}} [r] at {-10} 93
	\pinlabel {\mbox{$p-k$ cusps}} [l] at 105 93
	\pinlabel $\alpha$ at 238 18
	\pinlabel $\beta$ at 259 8
	\pinlabel $\delta_1$ at 241 56
	\pinlabel $\delta_2$ at 238 75
	\pinlabel $\delta_k$ at 273 117
	\pinlabel $\delta_{k+1}$ at 307 94
	\pinlabel $\delta_n$ at 275 38
	\endlabellist
\centering	
\includegraphics[scale=1]{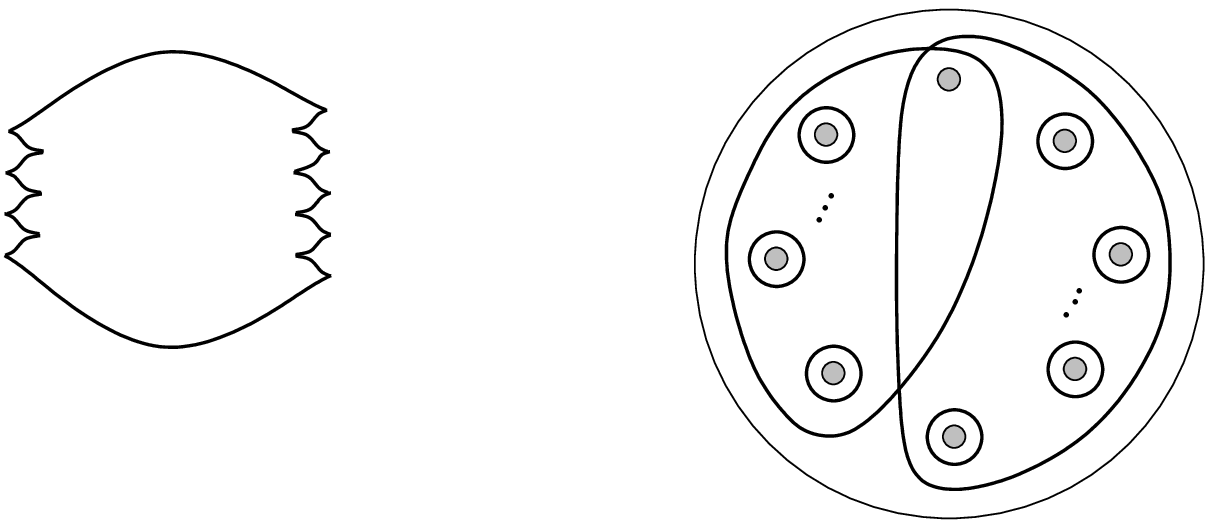}

\caption{An open book decomposition for the contact structure $\xi_k$ on $L(p,1)$ compatible with the surgery diagram on the left. 
(Here $p=n+1$.) The monodromy is the product $\Phi= D_{\alpha} D_{\delta_1} \dots D_{\delta_{k-1}} D_{\delta_{k+1}} D_{\delta_{n}} D_{\beta}$   of positive Dehn twists around the pictured curves.}
\label{lens-monodr} 
\end{figure}

In this setup, the page framing matches the Thurston--Bennequin 
framing of the Legendrian unknot; performing a positive Dehn twist on the curve representing the unknot, we obtain an open book decomposition 
corresponding to the result of Legendrian surgery. The resulting open book for $(L(p,1), \xi_k)$ is shown on Figure \ref{lens-monodr}: 
the page of the open book is a disk with $n=p-1$ holes, and the monodromy  
$$
\Phi= D_{\alpha} D_{\delta_1} \dots D_{\delta_{k-1}} D_{\delta_{k+1}} D_{\delta_{n}} D_{\beta} 
$$ 
is the product of positive Dehn twists $D_{\delta_i}$  around  each of the holes except the $k$-th, and 
the positive Dehn twists around the curves $\alpha$ and $\beta$. Here the order of Dehn twists is unimportant, 
because for a general open book it only matters up to cyclic permutation, and here the boundary  twists commute with all 
other Dehn twists. However, it is convenient to fix notation now:

\begin{convention} Throughout the paper, we adhere to braid notation for products of  Dehn twists: in the expression $D_{\alpha} D_{\beta}$,
$D_{\alpha}$ is performed first.  
\end{convention}

We first show that any positive factorization of $\Phi$ consists of Dehn twists enclosing the same holes.
(In a disk with holes, every simple closed curve separates, and we say that a collection of holes is enclosed 
by a curve if the holes lie in the component not containing the outer boundary of the disk. 
Abusing the language, we will often talk about holes enclosed by Dehn twists.)    

\begin{lemma} \label{factor} For the open book shown on Figure \ref{lens-monodr}, any positive factorization 
of the monodromy $\Phi$  must be given by the product of 
the Dehn twists $D_{\delta_1}, \dots, D_{\delta_{k-1}}, D_{\delta_{k+1}}, \dots D_{\delta_{n}}$,
and the Dehn twists $D_{\alpha'}$ and $D_{\beta'}$ around some curves $\alpha'$ and $\beta'$, such that $\alpha'$  encloses 
the same holes as $\alpha$, and $\beta'$ the same holes as $\beta$. 
\end{lemma}

To prove the lemma, we will need to look closely at the mapping class group of a planar surface; it will be convenient to work 
with its presentation given by Margalit and McCammond in \cite{Marg}. (Note:  The conventions in \cite{Marg} are opposite ours and those used by most 4-manifold or symplectic topologists; they use left-handed Dehn twists as the positive generators of the presentation.)

Consider a round disk with punctures arranged at the vertices of a regular $n$-gon contained in the disk; let the holes be small 
neighborhoods of the punctures, and denote by $\DD_n$ the resulting disk with $n$ holes.  We say that a simple closed curve in $\DD_n$ 
is convex if it is isotopic to the boundary of the convex hull of some of the holes; a Dehn twist around a convex curve is said to be convex. 
By \cite{Marg}, the mapping class group has a presentation 
with generators given by all (distinct) convex Dehn twists, and relations of the following two types. The first type states that Dehn twists around disjoint 
curves commute. The second type consists of all possible lantern relations; a lantern relation is the relation of the sort 
$D_A D_B D_C D_{A\cup B\cup C} = D_{A\cup B} D_{B \cup C} D_{A \cup C}$. Here  $A$, $B$, $C$ are disjoint collections of holes, 
$D_A$, $D_{A \cup B}$, etc, are convex Dehn twists around the curves enclosing the corresponding sets of holes, and the collections 
$A$, $B$, and $C$ are  
such that the cyclic clockwise ordering of all the holes in 
$A \cup B \cup C$, induced from their convex position on the disk, is compatible with the ordering where we list all holes from $A$ in their 
cyclic order, then all holes from $B$, then all holes from $C$.
Each collection $A$, $B$, $C$ may contain one or more holes; 
the ordering condition ensures that the Dehn twists are performed around curves arranged as in the usual lantern relation. See \cite{Marg} for details.

\begin{proof}[Proof of Lemma \ref{factor}] First, observe that although the set of holes enclosed by a simple closed $\gamma$ 
does not determine $\gamma$ up to isotopy, it determines, up to conjugacy, the class of the Dehn twist $D_{\gamma}$ in $\Map$. (If $\gamma$, 
$\gamma'$ enclose the same holes, the conjugacy between $D_{\gamma}$ and $D_{\gamma'}$ is given by the diffeomorphism $h\in \Map$ that 
takes $\gamma$ to $\gamma'$.) Consequently, the collection of holes enclosed by $\gamma$ uniquely determines the image of $D_{\gamma}$ in 
the homology $H_1 \Map$.  

It will be helpful to decompose all Dehn twists as follows: given a Dehn twist $D_{\gamma}$ around a convex curve $\gamma$ that encloses 
$r$ holes, we apply the lantern relation repeatedly to write $D_{\gamma}$ as a composition of Dehn twists around (convex) curves enclosing 
all possible pairs of holes among the given $r$, together with some Dehn twists around single holes. 
Each pair of holes will contribute exactly one positive 
convex Dehn twist into this decomposition; each hole will have $r-2$ negative Dehn twists around it. Similarly, we can decompose an arbitrary Dehn twist in the same fashion; however, the  Dehn twists enclosing pairs  of holes will no  longer have to be convex.  
We will refer to Dehn twists around a pair of holes as ``pairwise'' Dehn twists, and to Dehn twists around single holes as boundary twists.
(The Dehn twist around the outer boundary component of the disk {\em will not} be referred to as a boundary twist.)

Now we can decompose an arbitary element $\phi \in \Map$: write it as a product of Dehn twists, and decompose each of them as above. 
(We will not be recording the order of Dehn twists in the monodromy.) Let $m_{ij}(\phi)$ be the multiplicity of the positive Dehn twist containing only 
the $i$-th and $j$-th holes in the decomposition of $\phi$; similarly,  $m_i(\phi)$ be the multiplicity of the positive Dehn twist $\delta_i$ around 
the $i$-th hole. The collection of integers $m_{ij}, m_i$ is well-defined for $\phi \in \Map$: because the multiplicities  $m_{ij}, m_i$  
are invariant under lantern relations, they do not depend on the factorization of $\phi$. Moreover, the collection $\{m_{ij}, m_i\}$
uniquely determines the image of $\phi$ in  $H_1 \Map$. 
 
Let $\Phi$ be the monodromy of the open book for $(L(p,1), \xi_k)$ pictured in Figure \ref{lens-monodr}.  
Assume that $\Phi$ is factored into a product 
of some positive Dehn twists. Pick the $k$-th hole, which is enclosed by both curves $\alpha$ and $\beta$ in  Figure \ref{lens-monodr}, and 
consider the Dehn twists around it in the new factorization of $\Phi$. We consider only those Dehn twists that enclose at least two holes; 
suppose there are $l$ of them, and write  $n_1, n_2, \dots, n_l$ for the number of holes they enclose.

 Now, we compute the multiplicities of the pairwise and boundary twists for $\Phi$ that involve the $k$-th hole.
The Dehn twists $D_{\alpha}, D_{\beta}$ contribute  $n-1$ positive pairwise Dehn twists involving $k$-th hole   ($m_{i k}=1$ for every $i\neq k$) and 
$(k-2) +(n-k+1 -2)=n-3$ negative boundary twists around it.  On the other hand, the $l$ positive Dehn twists contribute $(n_1-1) + (n_2-1)+\dots +(n_l-1)$ 
positive pairwise Dehn twists involving the $k$-th hole, and $(n_1-2)+(n_2-2)+\dots (n_l-2)$ negative boundary Dehn twists. 
As the new factorization of $\Phi$
may have some positive boundary twists $D_{\delta_k}$  that we haven't yet taken into account, 
\begin{align}  
n-1 &= (n_1-1) + (n_2-1)+\dots +(n_l-1), \label{eq1}\\ 
n-3 & \geq (n_1-2)+(n_2-2)+\dots (n_l-2) \label{eq2}.
\end{align} 
 
Two cases are then possible: 

(i) $l=2$, and the new positive factorization of $\Phi$ has no boundary twists around the $k$-th hole, or     
 
(ii) $l=1$, and there is exactly one positive boundary twist $D_{\delta_k}$ around the $k$-th hole.
 
To rule out the second case, notice that since $m_{i k}(\Phi)=1$ for every $i\neq k$, the unique positive Dehn twist 
in the decomposition of $\Phi$ would have to enclose all $n$ holes. But then we would have $m_{1n}(\Phi)=1$, which contradicts 
the original definition of $\Phi$ from Figure~\ref{lens-monodr}.

The first case then tells us that the new positive decomposition of $\Phi$ has exactly  two twists enclosing more than one hole each;
there may also be some positive boundary twists. Examining the multiplicities $m_{ij}(\Phi)$, $m_i(\Phi)$ again, we see that one Dehn twist must be around a curve
that encloses the holes $1, 2, \dots k$, the other around a curve that encloses the holes $k, k+1, \dots n$;
denote the first curve by $\alpha'$, the second by $\beta'$.
 In addition, there must be one positive boundary twist $D_{\delta_i}$ 
around $i$-th hole for each $i\neq k$. \end{proof}

\begin{figure}[htb!] 
	\labellist
	\small \hair 2pt
	
	\pinlabel $0$ at 8 200
	\pinlabel $0$ at 177 267
	\pinlabel $0$ at 325 238
	\pinlabel $0$ at 315 25
	\pinlabel $0$ at 174 -10

	\pinlabel $-1$ at 119 194
	\pinlabel $-1$ at 251 255
	\pinlabel $-1$ at 359 188
	\pinlabel $-1$ at 346 60
	\pinlabel $-1$ at 218 9
	\pinlabel $-1$ at 118 120
	
	\pinlabel $\alpha'$ at 100 88
	\pinlabel $\beta'$ at 118 157

	\endlabellist
\centering
\includegraphics[scale=1.0]{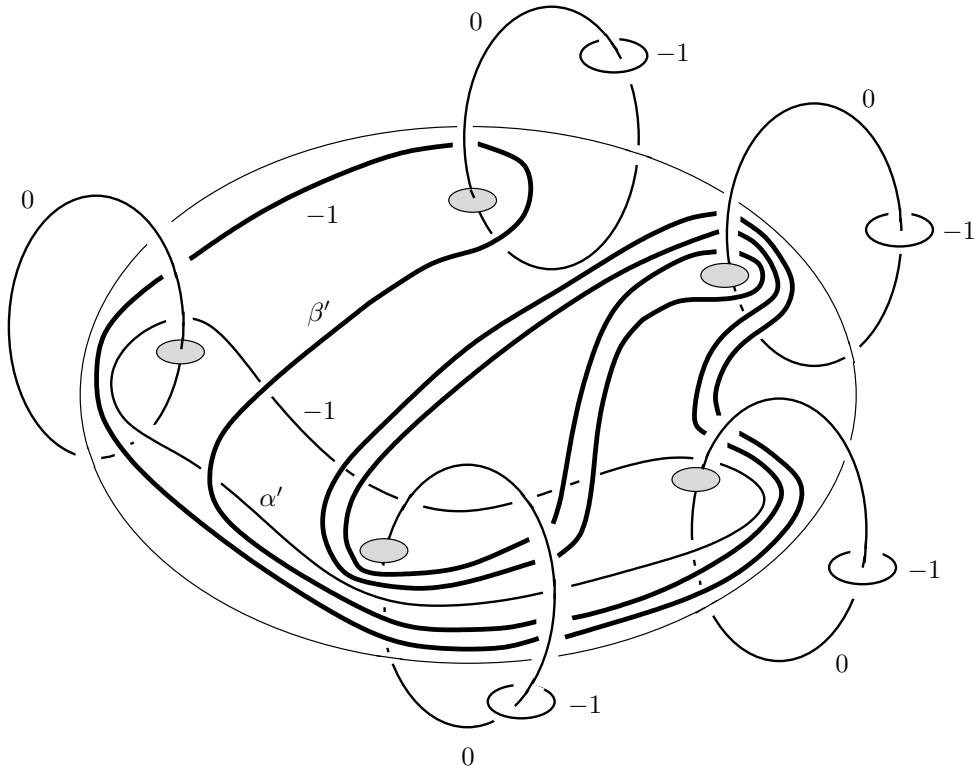}
\caption{Surgery diagram for the open book.}\label{obook} 
\end{figure}  

\begin{figure}[htb!] 
	\labellist
	\small \hair 2pt
	\pinlabel $0$ at 125 122
	\pinlabel $-1$ [l] at 146 78
	\pinlabel $+1$ [bl] at 285 116
	\pinlabel {\LARGE{-1}} at 396 101
	\endlabellist
\centering
\includegraphics[scale=1.0]{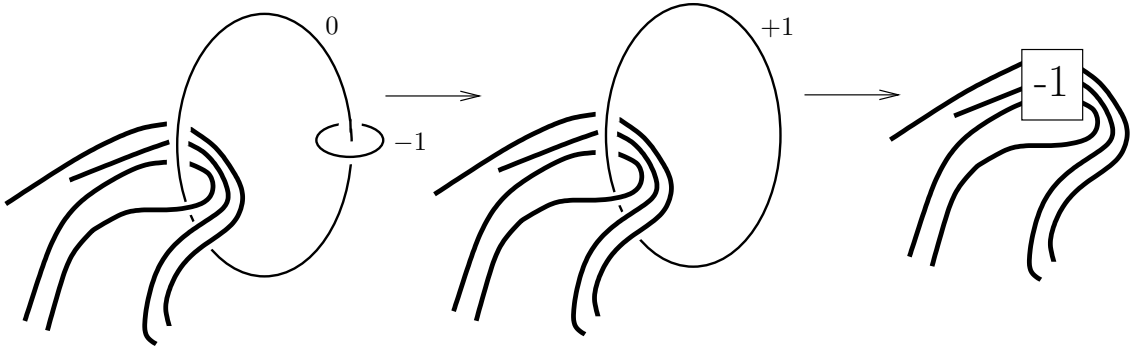}
\caption{Blowing down to destabilize at holes enclosed by $\beta'$.}\label{blowdowns} 
\end{figure}  
\begin{lemma} \label{frame} The open book whose monodromy is the product 
 $D_{\alpha'} D_{\delta_1} \dots D_{\delta_{k-1}} D_{\delta_{k+1}} \dots D_{\delta_{n}}$
represents $(S^3, \xi_{std})$. The knot in $S^3$ induced by $\beta'$ 
is the unknot; the framing on this unknot induced by the page framing of $\beta'$ is $-p+1$ (compared to the Seifert 
 framing).
\end{lemma}

\begin{proof} We can always find a self-diffeomorphism of the disk that maps $\alpha'$ to $\alpha$, 
so we may assume that the curve $\alpha'$ is standard. (Note, however, that we cannot simultaneously map 
$\alpha'$ to $\alpha$ and $\beta'$ to $\beta$). The open book with the monodromy  
$D_{\alpha'} D_{\delta_1} \dots D_{\delta_{k-1}} D_{\delta_{k+1}} \dots D_{\delta_{n}}$  
is obtained from the standard open book for $(S^3, \xi_{std})$ with annular pages by 
$n-1$ stabilizations, so it represents the standard tight $S^3$ as well.

\begin{figure}[htb!] 
	\labellist
	\small \hair 2pt
	
	\pinlabel 0 at 52 150 
	\pinlabel {-1} [r] at 0 131
	\pinlabel $k-p$ at 62 123
	\pinlabel {all +1} [bl] at -10 173
	\pinlabel {framed} [bl] at -10 166

	\pinlabel 0 at 204 156  
	\pinlabel -1 at 149 136
	\pinlabel +1 at 242 166
	\pinlabel +2 at 260 159
	\pinlabel +2 at 282 159 
	\pinlabel +2 at 302 159
	\pinlabel $k-p$ at 215 129

	\pinlabel $k-p$ [r] at 189 19
	\pinlabel +1 [b] at 220 38 
	\pinlabel +2 [b] at 242 38
	\pinlabel +2 [b] at 263 38
	\pinlabel +2 [b] at 285 38

	\pinlabel $-p$ at 64 7 
	
	\endlabellist
\centering
\includegraphics[scale=1.25]{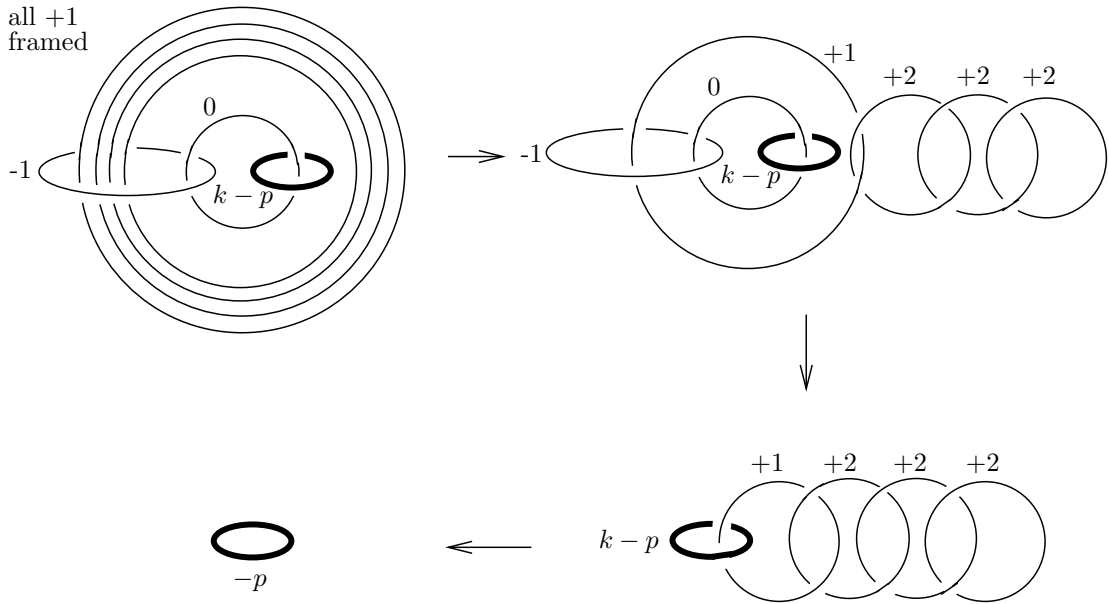}
\caption{Computing the framing. The thick circle represents the knot induced by $\beta'$, and is to be understood schematically: a priori it doesn't have to be an unknot, and its geometric linking with other components may be more complicated. 
However, the linking numbers  are as shown: the knot $\beta'$ has $lk=1$ with the 0-framed component, and zero linking numbers with the rest. It follows that the framings will change as dictated by Kirby moves 
in the picture.}\label{framing} 
\end{figure}  

 Next, we compute the page framing of the knot induced by $\beta'$; or rather, we will compute 
the surgery framing, which is the page framing minus 1. (We know that surgery on the corresponding knot yields 
$L(p,1)$, so the surgery framing must be $\pm p$, but the sign needs to be determined.)
Consider the Kirby diagram corresponding to the open book, Figure \ref{obook}.  We will 
destabilize the open book, starting with the holes enclosed by $\beta'$ and not by $\alpha'$. 
The knot  $\beta'$ will no longer lie on the page; 
in the Kirby diagram, destabilizations amount to blowdowns  shown on Figure \ref{blowdowns}.  After $n-k=p-1-k$ blowdowns, the framing 
decreases by $n-k$. The knot $\beta'$ has the linking number 1 with  
the $0$-framed unknot corresponding to the $k$-th hole (the one enclosed by both $\alpha$ and $\beta$), and 
the linking number zero with every other component of the surgery link. (Note that the geometric linking may be 
quite complicated.) Next, we perform further Kirby moves as in Figure \ref{framing}; knowing all 
linking numbers in the picture suffices for the framing calculation, even if the topological type of the knot $\beta'$ 
is unknown. It follows that the surgery framing of the original knot induced in $S^3$ by the curve  $\beta'$ is $-p$.
To see that $\beta'$ must be the inknot, we invoke a theorem of Kronheimer--Mrowka--Ozsv\'ath--Szab\'o \cite{KMOS} that states that the result of a $-p$-surgery on a knot in $S^3$ can be the lens space $L(p,1)$
only if the knot is the unknot.
\end{proof}

\begin{proof}[Proof of Theorem \ref{main}] Lemma \ref{frame}
implies that the Lefschetz fibration $X$  corresponding to the factorization of the monodromy as the product $D_{\alpha'} D_{\delta_1} \dots D_{\delta_{k-1}} D_{\delta_{k+1}} \dots D_{\delta_{n}} D_{\beta'}$  
is diffeomorphic to $D^4$  with a 2-handle attached along a $-p$ framed unknot. (Note that the order of Dehn twists in this product 
is not important.) This is the standard filling of $L(p,1)$, the disk bundle with 
Euler number $-p$. Moreover, $X$  has a Stein structure that arises from the Legendrian surgery on a Legendrian representative of the unknot 
$\beta'$, with $tb = -p+1$. Since Legendrian unknots are classified by their Thurston--Bennequin and rotation numbers \cite{EF}, we know 
that the only unknot that can produce $(L(p,1), \xi_k)$ is the one shown on Figure \ref{lens-monodr}, up to Legendrian isotopy.
Because a compatible symplectic structure on a Lefschetz fibration is unique up to symplectic deformation \cite{Go2}, it follows that 
all Stein structures on $X$ are sympectic deformation equivalent;  by Theorem \ref{wendl}, this means that the Stein  filling of 
$(L(p,1), \xi)$ is unique up to symplectic deformation, and any strong symplectic filling is unique up to symplectic deformation and blow-up. 

Finally, recall that every weak filling of a rational homology sphere can be modified into a strong filling \cite{OO}; it follows that the weak symplectic filling of  $(L(p,1), \xi)$ is also unique, up to deformation and blow-up. 
\end{proof}

\begin{remark} A very similar argument gives a new proof of McDuff's result \cite{McD} on uniqueness of filling for the standard (universally tight) 
contact structure on $L(p,1)$ for $p\neq 4$. The page of the corresponding open book is a disk with $n=p-1$ holes; the monodromy is the product 
of positive Dehn twists around the holes (one twist for each hole) and the positive Dehn twist around  the outer boundary component. Decomposing this monodromy, we get one positive Dehn twist around each pair of holes, and $n-3$ negative boundary twists for each hole. If there's a different 
positive factorization with $l$ non-boundary Dehn twists involving an arbitrary fixed hole and  enclosing respectively  $n_1, n_2, \dots n_l$ holes,
we decompose them as before to see that equations (\ref{eq1}) and (\ref{eq2}) must again hold. It follows that a positive factorization must either 
be the one we started with, or it must have two ``non-boundary'' Dehn twists involving each puncture, with no boundary twists. When $n\neq 3$,
the second case is not possible, because otherwise some ``pairwise'' Dehn twists would not be present in the decomposition. 
It follows that there is a unique factorization of the monodromy, with positive boundary twists around each hole and one positive twist that encloses 
all the holes.  (Unlike the case with two non-boundary twists considered above, no application of 
the deep result of \cite{KMOS} is needed here.) The classification of Legendrian unknots completes the proof as before; we see that the Stein filling is unique up to symplectic 
deformation. (This is weaker than Hind's result \cite{Hi}.)

\begin{figure}[htb] 
	\labellist
	\small \hair 2pt

	\pinlabel {$-1$} at 92 248
	\pinlabel {$-1$} at 124 197
	\pinlabel {$-1$} at 57 174

	\pinlabel {$-1$} at 334 174
	\pinlabel {$-1$} at 341 249
	\pinlabel {$-1$} at 396 199

	\pinlabel {$-2$} at 414 76
	\pinlabel {$-1$} at 409 27

	\pinlabel \Large{-2} at 232 42
	\pinlabel {$-2$} at 252 83
	\pinlabel {$-1$} at 252 26

	\pinlabel \Large{-2} at 56 41
	\pinlabel {$-3$} at 88 59

	\endlabellist
\centering
\includegraphics[scale=0.8]{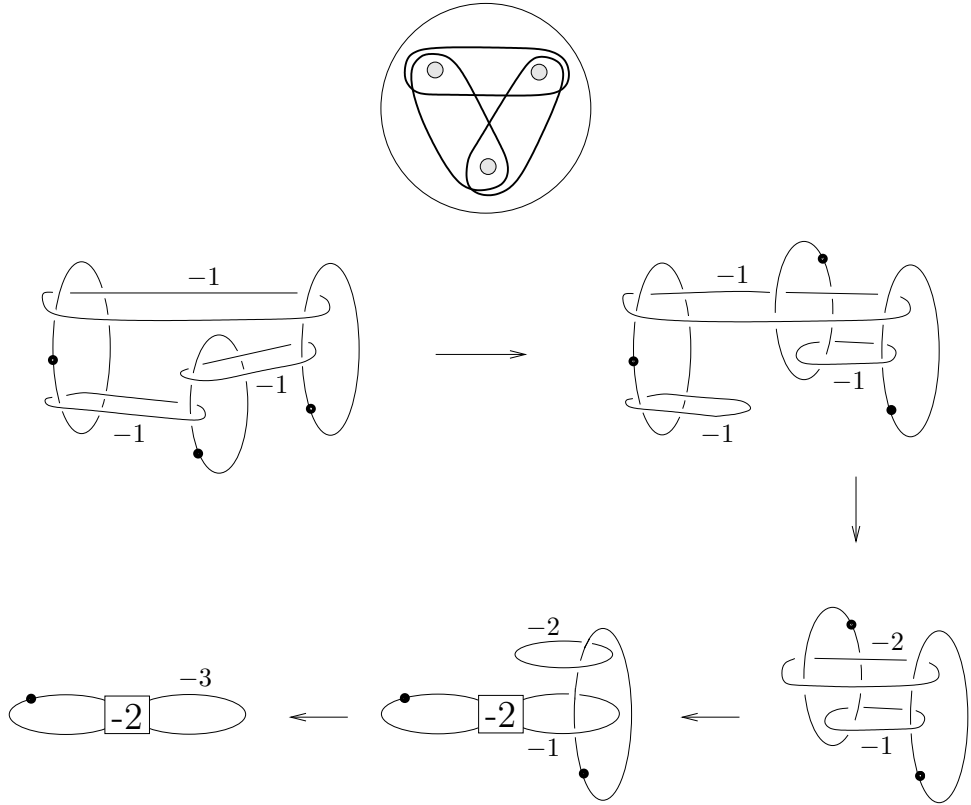}
\caption{Using the lantern relation to construct a non-standard Stein filling of $L(4,1)$.}\label{lantern} 
\end{figure}

When $n=3$ (i.e. $p=4$), 
an alternate positive factorization of the given monodromy indeed exists. It is given by the classical lantern relation, and corresponds to a non-standard filling of $(L(4,1), \xi_{std})$ which is a rational homology ball. (See Figure \ref{lantern} for a Kirby calculus picture
demonstrating that the lantern relation produces the non-standard filling constructed in \cite{McD}, \cite{Li}.) Note that we do not check that the non-standard 
filling is unique: indeed, the images of all the Dehn twists in $H_1 \Map$ are uniquely determined, but we do not have an appropriate analog of 
Lemma \ref{frame} for this case.   
\end{remark}

\begin{figure}[htb!] 
\includegraphics[scale=0.7]{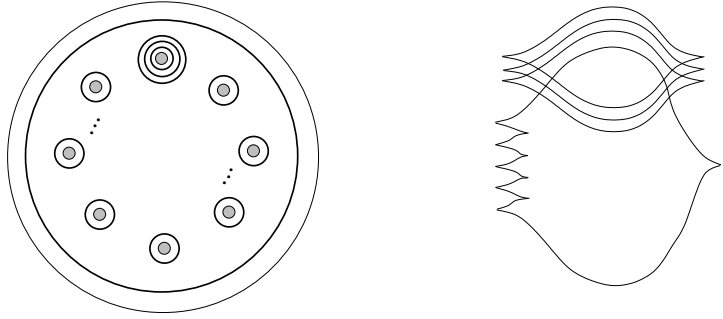}
\caption{Universally tight 
contact structures on $L(pk+1,p)$.}\label{more-lens} 
\end{figure}

\begin{remark} The same strategy proves  uniqueness of fillings for the universally tight 
contact structures on $L(pk+1,p)$, provided that $p, k\geq 1$, and either $k \leq p-2$, or $p=2$. See Figure~\ref{more-lens}, where  $p$ is the number of cusps in the surgery 
diagram (and the number of holes in the page of the open book); $k-1$ is the number of standard unknots 
(and the number of the multiple  boundary twists). If  $k>p-2>0$, a lantern-type relation (see Figure~\ref{genlantern}) can be used 
to construct a Stein  filling different from the one given by the surgery diagram. These contact structures are covered by Lisca's work \cite{Li}, but our
technique gives a slightly stronger result: uniqueness of filling up to symplectic deformation, not just diffeomorphism.  
\end{remark}

\section{Non-fillable contact structures on Seifert fibered spaces}

In this section we examine fillability of contact structures on the spaces $M(-1, r_1, r_2, r_3)$. All 
of these contact structures can be represented by planar open books; using Wendl's theorem, we can prove that a contact structure is non-fillable 
(in the strong symplectic sense) by showing that its monodromy admits no positive factorization in 
the mapping class group of a disk with holes. Since all of the spaces we consider are rational homology spheres, weak symplectic non-fillability 
follows as well.  

The following lemma will be used  repeatedly to control the number of Dehn twists in possible positive factorizations. 

\begin{lemma} \label{number-twists} Let $\phi \in \Map$ be given as a product of (positive or negative) Dehn twists. Suppose a hole $q$ is enclosed by $b$ boundary and $k$ non-boundary twists. Here $k=k_+ - k_-$, where $k_+$ (resp. $k_-$) is the number of positive
(resp. negative) non-boundary twists enclosing $q$; similarly, the sign of $b$ determines whether the boundary twists 
are positive or negative. 
Then in any positive factorization of $\phi$ there are no more than $k+b$ non-boundary Dehn twists enclosing $q$.   
\end{lemma}

\begin{proof} We generalize equations (\ref{eq1}), (\ref{eq2}). Suppose the $k_+$ positive Dehn twists around $q$ enclose $m_1, m_2, \dots m_{k_+}$ holes, while 
the $k_-$ negative Dehn twists enclose respectively $M_1, M_2, \dots M_{k_-}$ holes each. Assume that $\phi$ has a positive factorization where the hole $q$ is enclosed  by $l$ non-boundary positive Dehn twists around resp.  $n_1, n_2, \dots, n_l$ 
holes.   Computing the multiplicity
$m_q$ of the boundary twist around $q$ in the decomposition of $\phi$, as well as the pairwise multiplicities $m_{q q'}$ 
of all pairs that involve $q$, we have 
\begin{align*}  
& (m_1-1) +\dots +(m_{k_+}-1) -(M_1-1) -\dots -(M_{k_-}-1) = (n_1-1) +\dots +(n_l-1), \\ 
&      (m_1-2) +\dots +(m_{k_+}-2) -(M_1-2) -\dots -(M_{k_-}-2) -b    \leq (n_1-2)+\dots (n_l-2).
\end{align*} 
It follows that $k+b= k_{+}-k_{-}+b \geq l$. 
\end{proof}

\begin{example} \label{ex-Xi}
To avoid the more tedious analysis of cases that we'll need later on, we begin with an 
example taken from \cite[Figure 7]{GLS}, which is a tight contact structure $\Xi$ on the Seifert fibered space 
$M=M(-1; \frac12, \frac12, \frac1p)$.    The contact structure $\Xi$  is given by the surgery diagram on the left of Figure \ref{nonfill1};  
we can convert it into an open book on the right. (The order of Dehn twists is unimportant for this particular product.)    It is shown in \cite{GLS} that for $p>2$ $(M, \Xi)$   is not Stein fillable, 
and that it is not symplectically fillable for $p \not \equiv 2 \mod 8$. Since $(M, \Xi)$ can be represented by a planar open book, 
Stein non-fillability together with Wendl's work immediately implies that  
$(M, \Xi)$ is not symplectically fillable for all  $p>2$.  We now obtain an alternative quick proof of non-fillability, using our technique.

\begin{figure}[htb] 
	\labellist
	\small \hair 2pt

	\pinlabel $q_1$ at 281 128
	\pinlabel $q_2$ at 321 129
	\pinlabel $t$ at 301 136 
	\pinlabel $s_1$ at 243 85
	\pinlabel $s_{p-1}$ at 363 81
	
	\endlabellist
\centering	
\includegraphics[scale=1.0]{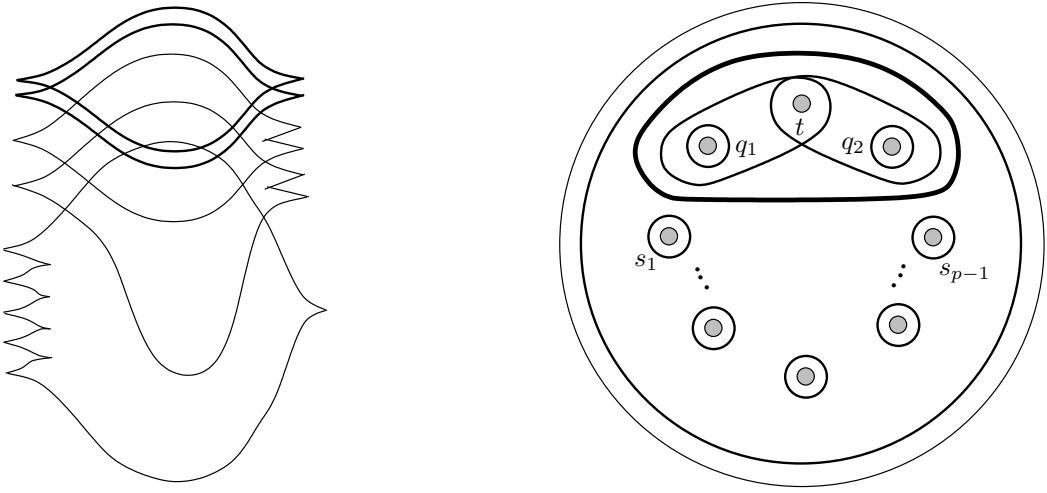}
\caption{The non-fillable tight contact structure on  $M(-1; \frac12, \frac12, \frac1p)$. The surgery diagram has two +1 contact surgeries
(on the thicker Legendrian unknots); the rest are Legendrian surgeries. Accordingly, the monodromy is the product of one negative Dehn twist
(around the thicker curve) and many positive ones.}\label{nonfill1} 
\end{figure} 

The monodromy $\Phi$ of the open book representing $(M, \Xi)$ is the product of several positive and one 
negative Dehn twist. There are $p-1\geq 2$ holes outside of the negative Dehn twist; denote them by $s_1, \dots, s_{p-1}$.
Suppose that $\Phi$ is factored into a product of positive Dehn twists around some curves; let 
$\{D_{\alpha}\}$ be the set of non-boundary twists in this factorization.    
Lemma \ref{number-twists} implies that every hole in the picture can be enclosed by no more than 
two of $D_{\alpha}$'s. Moreover, the pair of holes $q_1$ and $q_2$ is enclosed by one positive and one negative Dehn twist, so the multiplicity 
$m_{q_1 q_2}$ is zero. This means that none of $D_{\alpha}$'s can enclose both $q_1$ and $q_2$.

\begin{figure}[htb] 
	\labellist
	\small \hair 2pt
	
	\pinlabel $t$ at 97 136
	\pinlabel $s_1$ at 86 50
	\pinlabel $q_1$ at 52 97
	\pinlabel $q_2$ at 141 106
	\pinlabel $\alpha_1$ at 52 139
	\pinlabel $\alpha_2$ at 143 59
	\pinlabel $\alpha_3$ at 131 147

	\endlabellist
	\centering

\includegraphics[scale=1.0]{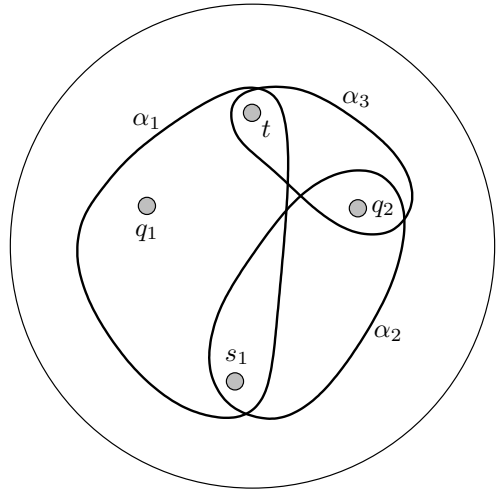}
\caption{Trying to construct a positive factorization for the monodromy of the open book from Figure \ref{nonfill1}.}\label{factorize} 
\end{figure} 
Now, consider the hole $s_1$. Since the pairs $\{q_1, s_1\}$,  $\{q_2, s_1\}$,  $\{t, s_1\}$  must all be enclosed with multiplicity 
1,  $q_1$ and $q_2$ cannot be enclosed together, and $s_1$ is enclosed by no more than two non-boundary twists, we must have 
a twist $D_{\alpha_1}$ that encloses  $t$, $q_1$, and $s_1$ (but not $q_2$), and another twist $D_{\alpha_2}$ that encloses  $q_2$ and $s_1$ (but not $q_1$ and $t$). 
(The roles of $q_1$ and $q_2$ may be reversed). See Figure \ref{factorize}. Because the pair $\{q_2, t\}$ has multiplicity 1, 
there must also be the twist $D_{\alpha_3}$ that encloses $t$ and $q_2$ but not $q_1$ and $s_1$. Next, consider the hole $s_2$.
Since $m_{s_1 s_2}=1$,  and  $D_{\alpha_1}$, $D_{\alpha_2}$
are the only non-boundary twists around $s_1$, the hole $s_2$ is enclosed by exactly one of $D_{\alpha_1}$ and $D_{\alpha_2}$.
If $s_2$ were in $D_{\alpha_1}$, we would have $m_{t s_2}= 2$ (if $s_2$ is  in $D_{\alpha_3}$),
or $m_{q_2 s_2}=0$ (if it isn't). Similarly, if $s_2$ were in $D_{\alpha_2}$, we have  $m_{t s_2}= 2$ or $m_{q_1 s_2}=0$. 
The contradiction shows that a positive factorization of $\Phi$ can't exist.  
\end{example}

\begin{remark} When $p=2$, the corresponding contact structure $\Xi$ is Stein fillable. This can be seen from the open book:
indeed, if there are no additional holes $s_2, \dots, s_k$, by the lantern relation Figure \ref{factorize} provides a positive factorization of the monodromy 
as the product $D_{\delta} D_{\alpha_3} D_{\alpha_2} D_{\alpha_1}$, where $D_{\delta}$ is the boundary twist around the hole $q_1$.
\end{remark}

The above argument readily generalizes to some more complicated  open books representing contact structures on other spaces $M(-1; r_1, r_2, r_3)$. 

\begin{figure}[htb!]
	\labellist
	\small \hair 2pt
	\pinlabel $s_1$ at 253 51
	\pinlabel $s_{p-1}$ at 345 49
	\pinlabel {$p-1$ twists} [tl] at 273 113
	\endlabellist
\centering
\includegraphics[scale=1.0]{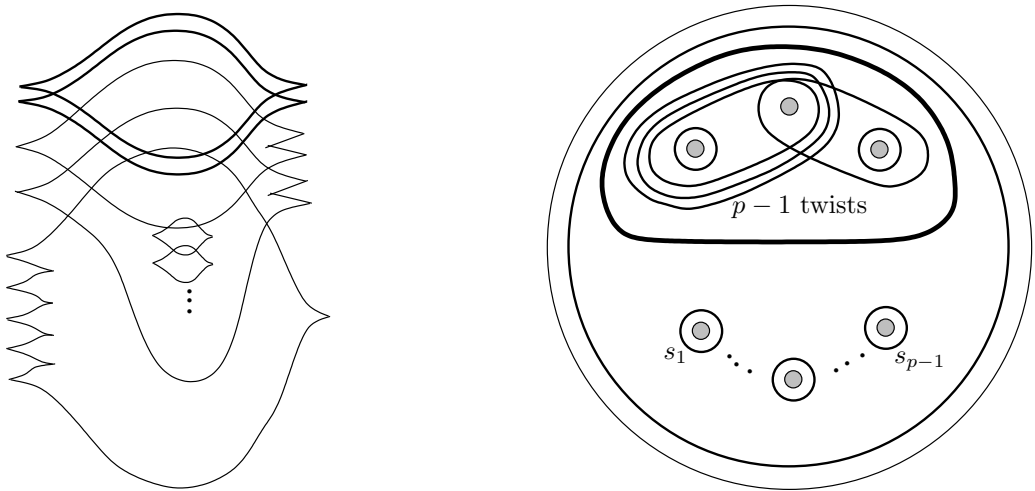}
\caption{A fillable contact structure on ${M(-1; \frac{p-1}{p}, \frac{1}{2}, \frac1{p})}$.}
\label{fillable} 
\end{figure}

 Recall that tight contact structures on the Seifert fibered space $M(-1; r_1, r_2, r_3)$, where $r_1, r_2\geq\frac12$,  
 were classified in \cite{GLS} and can be decribed as follows. The space  
 $M(-1; \frac12, \frac12, \frac1p)$ carries exactly three tight contact structures for each $p>2$; two of them are Stein fillable, 
 the third is the non-fillable contact structure $\Xi$ we studied in Example \ref{ex-Xi}. 
 It was shown in \cite{GLS} that  for  $r_1, r_2\geq\frac12$, all tight contact structures on
 $M(-1; r_1, r_2, r_3)$ can be obtained via Legendrian surgeries on these three 
 contact structures. The Legendrian surgeries arise from  
 the continued fraction expansions of $r_1, r_2, r_3$; thus, the surgeries are  performed on chains of Legendrian unknots, 
 with coefficients $a_i, b_i, c_i$. The Legendrian unknots  
 must be stabilized accordingly, to match $tb-1$ and the surgery coefficients. In general, many choices for stabilizations may be possible. 
 Fillability needs to be investigated only for those  contact structures obtained from $\Xi$, 
 as all others will automatically be fillable.

We first prove the ``fillable'' part of Theorem \ref{th2}. Indeed, suppose  
that $c_1=p \leq \max(k_1,k_2)$. We can assume that $k_1 \geq k_2$; this assumption implies that
$[a_1, a_2, \dots a_{p-1}]=[2,2,\dots,2]=\frac{p-1}{p}$. 
In that case, all the tight contact structures on $M(-1; r_1, r_2, r_3)$ 
are obtained from  $M(-1; \frac{p-1}{p}, \frac12, \frac1{p})$ by Legendrian surgeries. Thus, it suffices to prove 
the following 
\begin{prop} \label{fill} 
All tight contact structures on $M(-1; \frac{p-1}{p}, \frac12, \frac1{p})$ are Stein fillable.
\end{prop}

\begin{figure}[t!] 
	\labellist
	\small \hair 2pt

	\pinlabel {$\delta$} at 44 362
	\pinlabel {$\delta_1$} at 109 409
	\pinlabel {$\delta_2$} at 125 368
	\pinlabel {$\delta_{k+1}$} at 108 266
	\pinlabel {$k$ twists} at 45 313

	\pinlabel {$\alpha_1$} at 304 399
	\pinlabel {$\alpha_2$} at 313 375
	\pinlabel {$\alpha_{k+1}$} at 318 268
	\pinlabel {$\beta$} at 399 344

	\pinlabel {$\alpha$} at 279 188
	\pinlabel {$\alpha_1$} at 195 156
	\pinlabel {$\alpha_{m}$} at 230 74

	\pinlabel {$\beta_m$} at 272 99
	\pinlabel {$\gamma_m$} at 185 60 
	\pinlabel {$\delta_{k+1}$} at 218 28

	\endlabellist
	\centering
\includegraphics[scale=0.8]{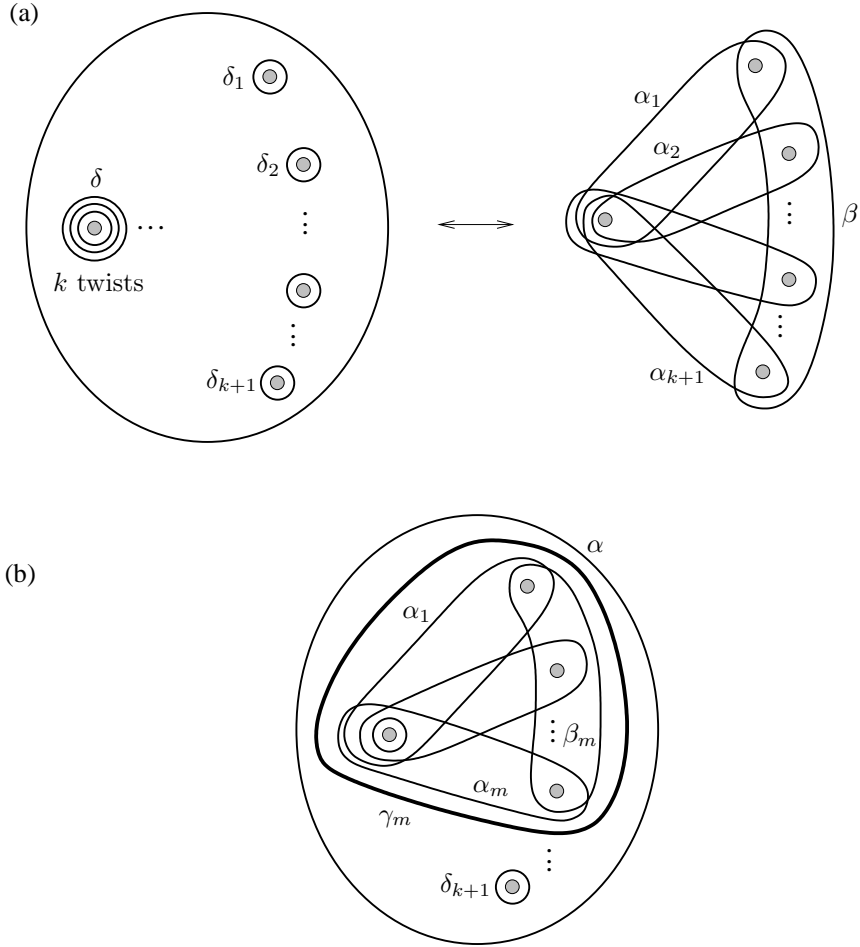}
\caption{(a) Generalized lantern relation for a disk with $k+2$ holes. (b) Proving it inductively: after $m-1$ applications of the classical lantern relation, 
the monodromy from the top right picture can be written as $k-m+1$ positive Dehn twists around $\delta$, a negative twist around $\gamma_m$, and 
positive twists around $\alpha_1,  \alpha_2, \dots, \alpha_{m}, \beta_m$, where the curve $\beta_m$ encloses $m$ holes.}
\label{genlantern} 
\end{figure}  

\begin{proof} We only need to check fillability of 
the (unique) contact structure on $M(-1; \frac{p-1}{p}, \frac{1}{2}, \frac1{p})$  obtained by Legendrian surgery on $\Xi$.
Its contact surgery diagram is shown on the left of Figure \ref{fillable}; to build the open book, we translate the chains of small 
unknots in the surgery diagram into sequences of push-offs  of the stabilized unknot they are linked to.  This is possible because the standard Legendrian meridian  and the push-off of a Legendrian knot are Legendrian isotopic in the manifold obtained by Legendrian surgery on this knot, \cite{DG}. The resulting open book is on the right of Figure \ref{fillable}.   
To prove the proposition, we rewrite the monodromy  as a product of positive 
Dehn twists using the lantern-type relation of Lemma \ref{lantrn} below. 
\end{proof}

\begin{lemma} \label{lantrn} In the mapping class group of the disk with $k+2$ holes, the relation 
$$
(D_\delta)^k D_{\delta_1} D_{\delta_2} \dots D_{\delta_{k+1}} D_{\alpha}=D_{\beta} D_{\alpha_{k+1}} \dots D_{\alpha_2} D_{\alpha_{1}}  
$$
holds for positive Dehn twists around curves shown on Figure \ref{genlantern}(a).  Our convention for products means that  
in the right-hand side, $D_{\beta}$ is performed first.   

\end{lemma} 
\begin{proof} The classical lantern relation, together with an inductive argument, shows that the relation  
$(D_\delta)^k D_{\alpha} D_{\delta_1} D_{\delta_2} \dots D_{\delta_k+1}  =  (D_{\delta})^{k-m+1} D_{\alpha} (D_{\gamma_m})^{-1} D_{\beta_m} D_{\alpha_m} \dots D_{\alpha_{2}} D_{\alpha_1}$ holds for each $m$, $k+1 \geq m \geq 2$. The right-hand side of 
this identity is illustrated on Figure \ref{genlantern}(b); note that $\gamma_{k+1}=\alpha$. 
\end{proof}

\begin{figure}[htb] 
	\labellist
	\small \hair 2pt

	\pinlabel {$k_1$ standard} [l] at -10 398
	\pinlabel {unstabilized} [l] at -10 390
	\pinlabel {unknots} [l] at -10 382

	\pinlabel {$k_2$ standard} [l] at 214 321
	\pinlabel {unstabilized} [l] at 214 313
	\pinlabel {unknots} [l] at 214 305

	\pinlabel {$n_1$ twists} at 30 198
	\pinlabel {$k_1$ twists} at 105 227
	\pinlabel {$k_2$ twists} at 210 225
	\pinlabel {$n_2$ twists} at 260 206 
	\pinlabel {$n_3$ twists} at 44 5
	
	\pinlabel $t$ at 143 166
	\pinlabel $q_1$ at 95 155
	\pinlabel $q_2$ at 191 155
	\pinlabel $w_1$ at 48 133
	\pinlabel $w_2$ at 242 132 
	\pinlabel $s_1$ at 101 80 
	\pinlabel $s_{c_1 - 1}$ at 195 80

	\endlabellist
	\centering

\includegraphics[scale=1.0]{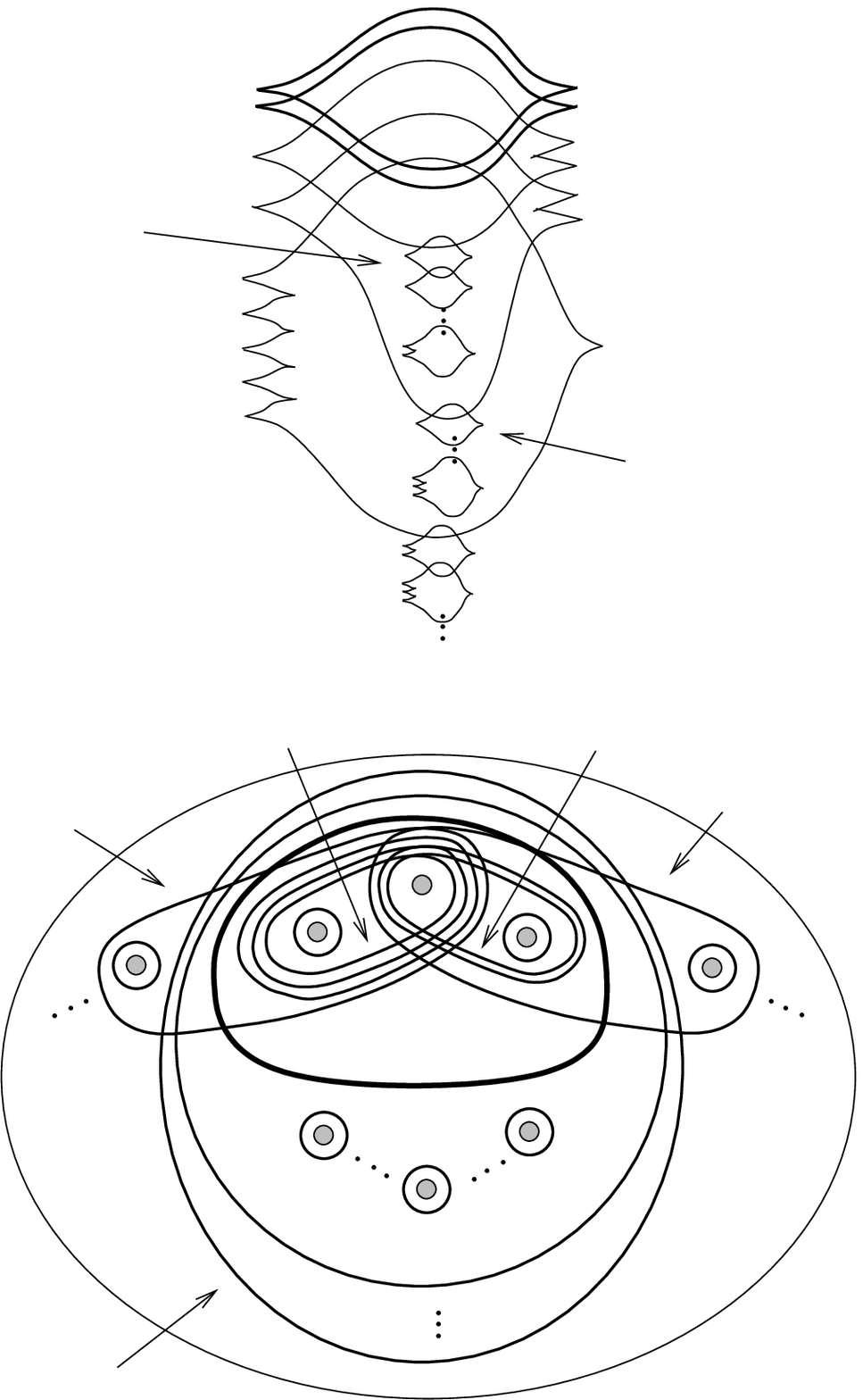}
\caption{Non-fillable tight contact structures on  $M(-1; r_1, r_2, r_3)$, with  $1> r_1, r_2 \geq \frac1{2}$.}\label{nonfill-new} 
\end{figure} 
It remains to prove the non-fillability part of Theorem \ref{th2}.
To this end, we construct the candidate non-fillable structures as follows. As before, we consider 
contact structures obtained by Legendrian surgery on $(M(-1; \frac12, \frac12, \frac1p), \Xi)$ (for various values  $p>2$). The coefficients $a_1, \dots, a_{k_1}, b_1, \dots b_{k_2}$ are equal to 2 and thus 
correspond to standard unknots with no stabilizations. 
If $k_1<n_1$ or $k_2<n_2$, there are some coefficients $a_{k_1+1}, b_{k_2+1}, \dots$ that are greater than 2. This means that the corresponding Legendrian unknots must be stabilized; we assume that all stabilizations are chosen {\em on the left}.  Similarly, we stabilize on the left all the unknots corresponding to the continued fraction expansion $r_3=[c_1, c_2, \dots c_{n_3}]$; note that 
$c_1=p$  if our contact structure is obtained from   $(M(-1; \frac12, \frac12, \frac1p), \Xi)$. The resulting contact structure $\xi$  on $M(-1; r_1, r_2, r_3)$  
is given by the contact surgery diagram and the open book shown on Figure \ref{nonfill-new}.
(As in Proposition \ref{fill},
we  use the identification of the meridian and the push-off of a Legendrian knot in a surgered manifold \cite{DG} to translate the chains of small unknots in the surgery diagram into sequences of iterated and possibly stabilized push-offs.)

Recall that $\xi$ is tight by \cite{GLS} (indeed, the Heegaard Floer contact invariant $c(\Xi)$ is non-zero, so  $c(\xi)\neq 0$ as well).

\begin{prop}  The contact structure $\xi$ on $M(-1; r_1, r_2, r_3)$ described above and shown on Figure \ref{nonfill-new} is non-fillable. \label{prop-nonfill}
\end{prop}

\begin{proof} We first point out several features of the open book representing $\xi$. The monodromy is the product of Dehn twists 
itemized below: there is a unique negative Dehn twist (explicitly mentioned in the list), and all other Dehn twists are positive.   
\begin{itemize}
\item There is a collection of $n_1$ Dehn twists enclosing the holes $q_1$ and $t$; $k_1$ of them enclose exactly these two holes,
the extra $n_1-k_1$, if present, all enclose a hole $w_1$, but not  $q_2$ or $s_1, s_2, \dots s_{c_1-1}$.
\item  There is a collection of $n_2$ Dehn twists enclosing the holes $q_2$ and $t$; $k_2$ of them enclose exactly these two holes,
an extra $n_2-k_2$, if present, all enclose  a hole $w_2$, but not $q_1$ or $s_1, s_2, \dots s_{c_1-1}$.
\item There is one negative Dehn twist enclosing three holes $q_1$, $q_2$ and $t$.
\item There are $n_3$ Dehn twists enclosing $c_1+2$ holes $q_1$, $q_2$, $t$, and $s_1, s_2, \dots s_{c_1-1}$, and perhaps 
some additional holes (but not $w_1$ or $w_2$).
\end{itemize}

Our task is to show that this monodromy admits no positive factorization. We will do so by analyzing the Dehn twists 
that can enclose the holes $s_i$. First, notice that since the multiplicity of any pair $\{w_1, s_i\}$,  $\{w_2, s_i\}$
is zero, any Dehn twist in a positive factorization that encloses of the $s_i$'s cannot enclose $w_1$ or $w_2$.
We will refer to the non-boundary Dehn twists that contain neither $w_1$ nor $w_2$ as {\em inner} Dehn twists; the (non-boundary) 
Dehn twists containing $w_1$ or $w_2$ will be called {\em outer} twists, and will not play much role in our argument.

First, we will calculate multiplicities of pairs of holes and  use Lemma \ref{number-twists} to establish the following properties 
of any possible positive factorization.

\begin{enumerate}[(a)]
\item exactly $n_3-1$ inner Dehn twists contain  both $q_1$ and $q_2$; all of these contain $t$
\item exactly $k_1$ inner Dehn twists contain $q_1$ and $t$, but not $q_2$ (call these $tq_1$-twists)
\item exactly $k_2$ inner Dehn twists contain $q_2$ and  $t$ but not $q_1$  (call these $tq_2$-twists)
\item there is at most one inner Dehn twist $D_{\alpha_1}$ that contains $q_1$ but not $q_2$ or $t$ 
\item there is at most one inner Dehn twist $D_{\alpha_2}$ that contains $q_2$ but not $q_1$ or $t$ 
\item each of the holes $s_1, s_2, \dots s_{c_1-1}$ is contained in exactly $n_3+1$ non-boundary Dehn twists, exactly $n_3$ of which enclose
$t$, exactly $n_3$ enclose $q_1$ and exactly $n_3$ enclose $q_2$  
\item  Any two holes $s_i$, $s_j$ are enclosed together by exactly $n_3$ Dehn twists
\end{enumerate}

To prove (a), observe that the multiplicity of the pair $\{q_1, q_2\}$ is $n_3-1$, which means there must be $n_3-1$ Dehn twists 
enclosing both holes. Each of them is an inner twist, since $m_{q_2, w_1}= m_{q_1, w_2}=0$. To see that all of these contain $t$, 
notice that by   Lemma \ref{number-twists} there are at most $n_1+n_2+n_3-1$ Dehn twists containing $t$; $n_1-k_1$ resp. $n_2-k_2$ 
of these must enclose $w_1$ resp. $w_2$, so there are at most $k_1+k_2+n_3-1$ inner Dehn twists containing $t$. Of these, at ledast $k_1+n_3-1$
contain $q_1$ (because $m_{q_1 t}=n_1+n_3-1$, and $q_1$ can be in no more than $n_1-k_1$ outer Dehn twists containing $t$ and $w_1$, 
and in none containing $w_2$). Similarly, at least $k_2$ contain $q_2$. Since all three holes $t$, $q_1$, and $q_2$ can be enclosed together by no more 
than  $n_3-1$ inner twists, the statements (a), (b) and (c) follow.

To prove (d), observe that by  Lemma \ref{number-twists} there are at most $n_1+n_3$ non-boundary Dehn twists enclosing $q_1$. Of these, 
$n_1-k_1$ are outer twists (containing $w_1$). This leaves at most $k_1+n_3$ inner twists enclosing $q_1$. By (a) and (b), at most one 
of these inner twists can contain neither $t$ nor $q_2$. The proof of (e) is similar. 
  
Finally, (f) follows from the fact that $m_{q_1 s_i}=m_{q_2 s_i}=m_{s_i t}=n_3$, together with (a) and Lemma \ref{number-twists};
(g) is merely a multiplicity count.   
  
We now show that properties (a)--(f) cannot hold  if $c_1-1>\max(k_1, k_2)$.  
  Indeed, (f), (a), (d) and (e) imply that each of the holes $s_i$ is enclosed by  $n_3-1$ twists containing the three holes $t$, $q_1$, $q_2$
(these are the twists described in (a)).
For the remaining two twists enclosing $s_i$, there are two possibilities: (i) a  $tq_1$-twist and $D_{\alpha_2}$, or (ii) a $tq_2$ twist and 
$D_{\alpha_1}$. Moreover, no two of the holes $s_i$ can be contained in the same $tq_1$-twist: they are also both contained 
in the $D_{\alpha_2}$ twist and in all the twists of (a), which contradicts (g). 
Similarly, no two $s_i$'s can be contained in the same $tq_2$-twist. Because $c_1-1>\max(k_1, k_2)$, we conclude that there must 
be  a non-empty subset of $\{s_1, s_2, \dots s_{c_1-1}\}$ for which (i) holds true, and a non-empty subset for which (ii) holds true.
Pick a hole from the first subset and another from the second; since we've listed all the $n_3+1$ twists enclosing each of these holes,     
we see that this pair of holes is enclosed together only  by the $n_3-1$ twists containing the three holes $t$, $q_1$, $q_2$. This is a contradiction with (g).
 \end{proof} 
    
\begin{remark}
In fact, we have established non-fillability for a much wider class of contact structures on Seifert spaces $M(-1; r_1, r_2, r_3)$ for which the condition 
$c_1-1>\max(k_1, k_2)$ is satisfied. Indeed, we only used the hypothesis that the unknots corresponding to $a_{k_1+1}$ and $b_{k_2+1}$ 
have at least one stabilization  
on the left, and that the unknot corresponding to $c_1$ has more than $\max(k_1, k_2)$ stabilizations on the left. All the other stabilizations 
may be chosen arbitrarily. 
\end{remark}

\section{Concluding remarks}

In the previous section, we worked with the spaces $M(-1, r_1, r_2, r_3)$ such that  $r_1, r_2\geq\frac12$
because this condition was essential for the classification results and for the proofs of tightness of \cite{GLS}.
In fact, with slightly more tedious case-by-case analysis, one can extend our non-fillability results 
to certain spaces $M(-1; r_1, r_2, r_3)$, with arbitrary $r_1, r_2 \in (0,1)$. We can also consider Seifert 
fibered spaces with more than 3 singular fibers; it is not hard to give sufficient conditions for an open book 
similar to the one in Figure \ref{nonfill-new} to represent a non-fillable contact structure.  One can also hope to understand 
fillability for a wider class of open books (not necessarily planar) by using the fiber sum construction of \cite{We2}. 
In many of these situations, contact structures can be shown to be non-fillable by easy combinatorial considerations in 
the homology of the mapping class groups of planar surfaces.

However, non-fillability does not seem to be an interesting property unless the contact structure is known to be tight.
In most cases, a proof of tightness requires an application of Heegaard Floer homology. For the spaces $M(-1, r_1, r_2, \dots, r_k)$, 
non-vanishing of the Heegaard Floer contact invariant (and thus tightness) can be checked via the criterion of 
\cite[Theorem 1.2]{LS}, although this often requires lengthy calculations related to Heegaard Floer homology (see \cite{LS2}).

It would be interesting to find hypotheses on the monodromy of a planar open book that ensure tightness of the corresponding contact 
structure. One related condition is the right-veering property of an open book \cite{HKM1}; indeed, a contact structure is tight 
if and only if every compatible open book is right-veering. In general, one needs to consider {\em all} compatible open books;
indeed, by stabilizing an arbitary open book, one can always obtain a  right-veering open book representing the given contact 
structure. This is shown in  \cite[Proposition 6.1]{HKM1}; in fact, the proof of that proposition shows that stabilizations can 
be done without increasing the genus of the open book (but increasing the number of its boundary components). Thus,  every contact structure 
supported by a planar open book can be supported by a right-veering planar open book. In fact, one can perform additional stabilizations
to increase the pairwise multiplicities of the monodromy; it is not hard to show that every contact structure supported by a planar open 
book is supported by a planar open book with right-veering monodromy and positive pairwise multiplicities. By contrast,  recall  
that a contact structure supported by an open book of genus one with one boundary component is tight if this particular open book is 
right-veering \cite{HKM2}; one can hope that since planar open books are another very special case, some sufficient conditions for 
a {\em given} monodromy to ensure tightness can be found.

In another direction, it would be interesting to generalize the classification results of Section 1. Indeed, it is possible to analyze 
positive factorizations in the abelianization of the mapping class group for a wide class of contact structures. However, we have no 
analogs of Lemma \ref{frame} (in fact no analogs of the theorem of \cite{KMOS}), and this poses a major obstacle for further classification results.

 \end{document}